\documentclass[a4paper]{article} 
\usepackage[cmfonts]{iba} 
\usepackage{booktabs} 
\usepackage{amsmath} 
\usepackage{amssymb} 
\usepackage{enumerate} 
\usepackage{graphicx,color} 
\usepackage{pgf}
\theoremstyle{plain}
\newtheorem{claim}{Claim}
\newcommand{\B}{\{0,1\}} 
\DeclareMathOperator\conv{conv} 
 
\DeclareMathOperator\Proj{Proj}

\def\dotcup{\mathbin{\setbox0\hbox{$\cup$}\rlap{\copy0}\raise.3\ht0\hbox
    to\wd0{\hfil$\cdot$\hfil}}} 

\usepackage{enumerate} 
 
\let\saveexample=\example
\renewcommand\example{\par\bgroup\small\saveexample}
\let\saveendexample=\endexample
\renewcommand\endexample{\saveendexample\egroup}

\usepackage[breaklinks,colorlinks,citecolor=blue]{hyperref} 
 
\renewcommand{\O}{O} 
 
\title{Intermediate integer programming representations using value disjunctions\footnote{Address: Otto-von-Guericke-Universit\"at  
    Magdeburg, Department of Mathematics/IMO, Uni\-ver\-si\-t\"ats\-platz~2,  
    39106 Magdeburg, Germany.\newline E-mail addresses: 
    \texttt{\{mkoeppe,\,louveaux,\,weismant\}@imo.\penalty0  
      math.\penalty0 uni-\penalty0 magdeburg.de}\newline
  The authors gratefully acknowledge support from the European TMR network ADONET 504438. 
}}
\author{Matthias K\"oppe\and Quentin Louveaux \and Robert Weismantel} 
\date{$\relax$Date: 2005/11/14 14:36:00 ${}$ -- ${}$Revision: 3.141 $ $} 
\begin{document} 
\maketitle 

\begin{abstract} 
  We introduce a general technique to create an extended formulation of a
  mixed-integer program.  We classify the integer variables into
  blocks, each of which generates a finite set of vector values.  
  The extended formulation is constructed by creating a new binary
  variable for each generated value.  Initial experiments show that the
  extended formulation can have a more compact complete description than
  the original formulation.  

  We prove that, using this reformulation technique, the facet description
  decomposes into one ``linking 
  polyhedron'' per block and the ``aggregated polyhedron''.  Each of these
  polyhedra can be analyzed separately.  For the case of identical
  coefficients in a block, we provide a complete description of the linking
  polyhedron and a polynomial-time separation algorithm.  Applied to the
  knapsack with a fixed number of distinct coefficients, this theorem
  provides a complete description in an extended space with a polynomial number
  of variables.

  Based on this theory, we propose a new branching scheme that analyzes the
  problem structure.  It is designed to be applied in
  those subproblems of hard integer programs where LP-based 
  techniques do not provide good branching decisions.
  Preliminary computational experiments show that it is successful for some
  benchmark problems of multi-knapsack type. 
\end{abstract} 
 
\section{Introduction} 

Extreme representations of the feasible points of a mixed-integer linear
optimization problem
 are either given by means of the facet defining inequalities in the
 original space or by a set of feasible mixed integer points whose convex hull contains the feasible region. It is well known that in principle one such extreme representation can be transformed into  the other extreme representation.
However from an algorithmic point of view both extreme representations are
very hard to achieve. 

This suggests to search for other, ``intermediate'' representations that are algorithmically more tractable, in the sense that they 
\begin{itemize}
\item require  less variables than the extreme representation by the vertices,
\item require  less constraints compared to the total number of facets of the convex hull,
\item have a simpler combinatorial constraint structure  than the facets of the convex hull in the original space and hence, the separation problem in the extended space is easier to solve. 
\end{itemize}
Intermediate representations of the feasible region are complete descriptions 
of an extended formulation of the original problem.  To make this notion
precise, we define:
\begin{definition}[Representation by projection]
Let $P\subseteq \R^n$, $P' \subseteq \R^d$ be two rational polyhedra and $B\in \Q^{n \times d}$ a rational matrix.
We call $P'\cap \Z^d$ a \emph{representation} of $P\cap \Z^n$ if the following two properties hold:
\begin{itemize}
\item[(a)] $P\cap \Z^n = \{\,\ve x \in \Z^n : \ve x=B\ve y,\; \ve y \in P' \cap \Z^d\,\}$. 
\item[(b)] $\conv(P\cap \Z^n) = \{\, \ve x \in \R^n : \ve x=B\ve y,\; \ve y \in P'\,\}$. 
\end{itemize}
Such a representation is called \emph{extreme} if either $d=n $ and $B=I$ or if 
$P' = \{\,\ve y \in \R_+^d : \sum_{i=1}^d y_i = 1\,\}$; otherwise, it is
called \emph{intermediate}.
\end{definition}
We remark that R.~K.~Martin \cite{rkmartin-1987} calls the sets $P'\cap\Z^d$
and $P\cap \Z^n$ ``strongly equivalent'' in this situation. 

In the literature, there are a couple of interesting examples of this type.
Chopra and Rao \cite{ChopraRao94,ChopraRao94b} introduced a directed formulation for the Steiner tree problem and showed that exponentially many inequalities in the undirected formulation are projections of a small number of directed inequalities.  
R.~K.~Martin~\cite{rkmartin:91} reports on the minimum spanning tree problem, 
which has as an inequality formulation of size $\O(2^n)$. It can, however, alternatively be described as the 
projection of an extended formulation which requires $\O(n^3)$ variables and $\O(n^2)$ 
constraints. 
Moreover, there are many further compact extended formulations for specific
combinatorial optimization problems, in particular for lot-sizing and
fixed-charge network problems; see, for instance, \cite{krarup-bilde-1977,
  rardin-choe-1979, rkmartin-1987}.

Next we illustrate on an example that also quite general problems such as knapsack problems can sometimes be described in an extended space such that 
the higher dimensional polyhedron is much more appealing than the original facet description.

\begin{example}
Consider the set of $\ve x\in\{0,1\}^8$ such that
  \begin{equation}
\label{knapex1}
    8x_0 - x_{1} -2x_{2} - 3x_{3} -4x_{4} -5x_{5} -6x_{6} -7x_{7} {}\leq0.
  \end{equation}
The convex hull of solutions to this knapsack problems is given by the following system of thirteen inequalities:
{\small
\begin{displaymath}
  \begin{array}{*9{@{}r}}
 x_0                    &        &          & {}- x_{3}                   &                             & {}- x_{5}                    & {}- x_{6}                     & {}- x_{7}                     & {}\leq 0 \\
 x_0                    &        &          &                             & {}- x_{4}                   & {}- x_{5}                    & {}- x_{6}                     & {}- x_{7}                     & {}\leq 0 \\
 x_0                    & {}-x_{1} & {}-x_{2} &                             &                             & {}- x_{5}                    & {}- x_{6}                     & {}- x_{7}                     & {}\leq 0 \\
 x_0                    & {}-x_{1} &          & {}- x_{3}                   & {}- x_{4}                   &                              & {}- x_{6}                     & {}- x_{7}                     & {}\leq 0 \\
 x_0                    &        & {}-x_{2} & {}- x_{3}                   & {}- x_{4}                   & {}- x_{5}                    &                               & {}- x_{7}                     & {}\leq 0 \\
 x_0                    &        & {}-x_{2} & {}- x_{3}                   & {}- x_{4}                   &                              & {}- x_{6}                     & {}- x_{7}                     & {}\leq 0 \\
 x_0                    & {}-x_{1} & {}-x_{2} & {}- x_{3}                   & {}- x_{4}                   & {}- x_{5}                    & {}- x_{6}                     &                               & {}\leq 0 \\
  {\bf2}x_0  & {}-x_{1} & {}-x_{2} & {}- x_{3}                   & {}- x_{4}                   & {}- x_{5}                    & {}- x_{6}                     & {}- x_{7}                     & {}\leq 0 \\
  {\bf2}x_0  &        & {}-x_{2} & {}- x_{3}                   & {}- x_{4}                   & {}- x_{5}                    & {}- x_{6}                     & {}-  {\bf2}x_{7}   & {}\leq 0 \\
  {\bf2}x_0  & {}-x_{1} &          & {}- x_{3}                   & {}- x_{4}                   & {}- x_{5}                    & {}-  {\bf2}x_{6}   & {}-  {\bf2}x_{7}   & {}\leq 0 \\
   {\bf3}x_0 & {}-x_{1} & {}-x_{2} & {}- x_{3}                   & {}- x_{4}                   & {}-  {\bf2}x_{5}  & {}-  {\bf2}x_{6}   & {}-  {\bf2}x_{7}   & {}\leq 0 \\
   {\bf3}x_0 & {}-x_{1} & {}-x_{2} & {}-  {\bf2}x_{3} & {}-  {\bf2}x_{4} & {}- x_{5}                    & {}-  {\bf2}x_{6}   & {}-  {\bf2}x_{7}   & {}\leq 0 \\
\   {\bf 5}x_0   & {}-x_{1} & {}-x_{2} & {}-  {\bf2}x_{3} & {}-  {\bf2}x_{4} & {}-   {\bf3}x_{5} & {}- {\bf4}x_{6} & {}- {\bf4}x_{7} & {}\leq 0
\end{array}
\end{displaymath}}%
One way to obtain an extended formulation for (\ref{knapex1}) is to introduce two new variables for the subsets $\{1,2\}$ and
$\{3,4\}$. This requires to introduce two new variables $x_{\{1,2\}}$ and $x_{\{3,4\}}$ which are equal to one if both elements $1$ and $2$ ($3$ and $4$, respectively) are selected.
This yields the following reformulation:
    \begin{displaymath} 
      \begin{array}{*{11}{@{}r}}
        8x_0 &{}- x_{1} &{}-2x_{2} &{}- 3x_{3} &{}-4x_{4} &{}-5x_{5} &{}-6x_{6} &{}-7x_{7} &{}-3x_{\{1,2\}} &{}-7x_{\{3,4\}} &{}\leq0\\
        &x_{1}&{}+x_{2}&&&&&&{}+  {   x_{\{1,2\}}}&&{}\leq1\\
        & & &x_{3}&{}+x_{4}&&&&&{}+\   {   x_{\{3,4\}}}&{}\leq1
      \end{array}
    \end{displaymath}
The convex hull of all feasible binary solutions to this system is given by the following list of nine inequalities: 
    {\small 
      \begin{displaymath}
        \begin{array}{*{11}{@{}r}}
          x_0    &        & &        & &{}-x_{5}& {}-x_{6} & {}- x_{7}& {} &{}-    x_{\{3,4\}}& {}\leq 0                                     \\
          x_0    &  {}-x_{1}& {}-x_{2} & &        & {}-x_{5} & {}-x_{6} & {}- x_{7}  & {}-   x_{\{1,2\}}  & & {}\leq 0                       \\
          x_0    &       & &{}-x_{3}& {}-x_{4} & &{}-x_{6}& {}- x_{7}&  {}-   x_{\{1,2\}}   & {}-    x_{\{3,4\}}& {}\leq 0                         \\
          x_0    &       & {}-x_{2} & {}-x_{3} & {}-x_{4} & {}-x_{5} & &{}- x_{7} &  {}-   x_{\{1,2\}}  & {}-    x_{\{3,4\}}& {}\leq 0             \\
          x_0   & {}-x_{1}& {}-x_{2} & {}-x_{3} & {}-x_{4} & {}-x_{5} & {}-x_{6}   & {} & {}-   x_{\{1,2\}}  & {}-    x_{\{3,4\}}& {}\leq 0             \\
            {\bf2}x_0   & {}-x_{1}& {}-x_{2} & {}-x_{3} & {}-x_{4} & {}-x_{5} & {}-x_{6}   & {}- x_{7}  & {}-   x_{\{1,2\}}  & {}-    x_{\{3,4\}}& {}\leq 0 \\
            {\bf2}x_0   &       & {}-x_{2} & {}-x_{3} & {}-x_{4} & {}-x_{5} & {}-x_{6}   & {}-  {\bf2}x_{7}  &  {}-   x_{\{1,2\}}  & {}-  {\bf2}   x_{\{3,4\}}& {}\leq 0 \\
          &       & &{}+x_{3}& {}+x_{4} & &        &        & &{}+    x_{\{3,4\}}& {}\leq 1                                                  \\
          & x_{1}& {}+x_{2} & &        & &        &        & {}+   x_{\{1,2\}}   & & {}\leq 1
        \end{array}
      \end{displaymath}}%
Note that not only the number of inequalities for the extended  formulation is
smaller than in the original space. More importantly, the structure of the
inequalities in the extended space is significantly nicer when compared to the
structure of the inequalities in the original space. For instance, the maximum
coefficient occuring in the inequalities in the higher dimensional space
is~$2$, whereas the highest coefficient in the inequalities in the original
space is already~$5$.   
\end{example} 

In the example, the extended formulation that we propose is based on introducing two new variables that correspond to products of variables in the original space. This is a special case of an extended formulation that one can obtain from the so-called Lift-and-Project approach.  This approach has its roots in the work of Egon Balas on disjunctive optimization  \cite{Balas71,Balas75a}. 
It was further refined in \cite{SheraliAdams90, LovaszSchrijver91,
  BalasCeriaCornuejols93, BienstockZuckerberg02} by introducing hierarchies of extended formulations whose variables represent more general subsets of original variables. 
The disadvantage of this approach is that the number of variables grows exponentially with the size of the subsets for which we introduce new variables.

The tool that we propose in this paper to generate an extended formulation is
the {\em value-disjunction procedure}. It is another generalization of
introducing  one  new variable for each  pair of original variables. However,
it  also applies to subsets of original variables of larger cardinality and offers lots of freedom in generating the extended formulation. It is a general way to 
produce intermediate representations for mixed integer optimization problems. In fact, it provides a hierarchy of new formulations. 
Specifically, for any subset of original variables we can always introduce an extended formulation that keeps the number of new variables linear in the size of the subset.

We introduce the value disjunction procedure in Section 2.  We then describe the convex hull of  
the given mixed-integer set as the intersection of several simpler polyhedra using the variables 
of the extended space. This is the \emph{structure theorem} for the value disjunction procedure. 
In Section~3 we  introduce the family of linking polyhedra. In the special but
important  case that such a linking polyhedron comes from the unweighted sum
of a set of variables, 
we completely describe the polyhedron by means of linear inequalities and equations. 
As an application of the structure theorem in Section~2 together with the polyhedral
characterizations of Section~3, we are able to determine an explicit description
of the convex hull of all solutions to a 0/1 knapsack problem with only a
fixed number of different weights. This is the topic of Section 4. 

Finally, in Section~5, we investigate one way of making computational use of
value disjunctions: By branching also on the new binary variables of the
extended formulation instead of only on the original variables, it is possible
to take more flexible branching decisions.  In fact, we propose such a
branching scheme for situations where none of the usual LP-based variable
selection criteria provides a solid basis for taking a branching decision.
Such situations frequently occur in very hard integer programs like the
market-split instances~\cite{corn-dawande:98}.  We investigate the effect of
branching simplifying the facet description: A branching decision is
considered good if the facet descriptions of the generated subproblems are 
significantly simplier than the original facet description.
Using experiments with randomly generated problem instances, we show that it is possible to make a branching
decision based on the structure of the problem which is better than
branching on the original variables.  
Finally we report on simple computational experiments
with a few hard integer programs, where we branch explicitly on the new binary
variables and then solve the subproblems with the branch-and-cut system
{\small CPLEX}.  
We obtain a significant reduction in both the number of nodes and the computation time.

\section{Value disjunctions} 

In this section, we present a structural result about an extended formulation
of a given mixed integer programming model. To this end, consider a bounded mixed-integer set of the form 
\begin{equation*} 
\mathcal F = \biggl\{\,(\ve x,\ve w) \in \Z^{n}_+\times \R^{d}_+ :
\sum_{j=1}^{n} A_j x_j + \sum_{j=1}^{d} G_j w_j \leq \ve b,\ \ve x\leq\ve u\,\biggr\}, 
\end{equation*} 
where $A_j$, $G_j\in \R^m$ for all $j$, $\ve b \in \R^m$, and $\ve u\in\Z_+^n$. We set $P=\conv \mathcal F$.  

Let us partition the set $N=\{1,\ldots, n\}$ into subsets $N_1, \ldots, N_K$. 
For each of the sets $N_i$, we determine all the possible vectors (``values'') generated
by the columns~$A_j$ belonging to the variables indexed by $N_i$:
$${\cal A}_i = \biggl\{\,\sum_{j \in N_i} A_j x_j : x_j \in \{0,\dots,u_j\}\text{ for
  $j\in N_i$}\,\biggr\}.$$
Since the integer variables are assumed to be bounded, the set ${\cal A}_i$
is finite; its cardinality $n_i = |{\cal A}_i|$ is at most $\prod_{j\in N_i} (1+u_j)$.  Let the
elements of ${\cal A}_i$ be numbered, ${\cal A}_i=\{\ve f^{N_i}_1,\ldots,\ve f^{N_i}_{n_i}\}$.
We shall associate with $\ve f^{N_i}_k$ a new binary variable 
$y^{N_i}_k$. In order to simplify the subsequent expositions, we shall also use the 
abbreviating notations $A(\ve x^{N_i}) = \sum_{j\in N_i} A_j$, and moreover $A(\ve
y^{N_i}) = \sum_{k=1}^{n_i} y^{N_i}_k \ve f^{N_i}_k$ and $A(\ve y) = \sum_{i=1}^K A(\ve y^{N_i})$.

We come to two major definitions that we make use of in this paper.  
\begin{definition} 
For a given subset $N_i$, we define the \emph{linking polyhedron} as
\begin{align} 
V_i=\conv \biggl\{\, (\ve x^{N_i},\ve y^{N_i}) \in \Z_+^{|N_i|}\times \B^{n_i}:& 
\sum_{j\in N_i} A_j x_j = \sum_{k=1}^{n_i} \ve f^{N_i}_k y^{N_i}_k\notag\\ 
& \sum_{k=1}^{n_i} y^{N_i}_k = 1 \notag\\  
& 0\leq x_i\leq u_i, \;i=1,\ldots, n \qquad  \biggr\}.\label{set Vi} 
\end{align} 
Furthermore we define the \emph{aggregated polyhedron} as
\begin{align} 
Q = \conv\biggl\{\,(\ve y,\ve w)&\in\B^{n_1 + \cdots + n_K}\times \R^d_+:\notag\\
&\sum_{i=1}^K \sum_{k=1}^{n_i} \ve f^{N_i}_k y^{N_i}_k + \sum_{j=1}^d G_jw_j\leq\ve b \notag\\ 
& \sum_{k=1}^{n_i} y^{N_i}_k \leq 1 \; \text{for all }i=1,\ldots, K \qquad\quad
\biggr\}.\label{set Q} 
\end{align} 
\end{definition} 

Thus, for every value $\ve f^{N_i}_k$ in a set $\mathcal A_i$ we are introducing  a
new binary variable $y^{N_i}_k$.  With this family of new variables, we can 
obtain a new, extended formulation of $\mathcal F$ by linking  
the original variables $\smash{x_j}$ with the new ``value variables'' $y^{N_i}_k$.
The precise link between the extended formulation and the original formulation
is given in the following theorem.  Before stating the theorem we
illustrate our constructions on an example. 
\begin{example}
Consider the convex hull $P$ of all binary solutions to the inequality 
\begin{displaymath}
  3x_{1} +3x_{2} + 3x_{3} +3x_{4} +4x_{5} +7x_{6} +8x_{7} +9x_{8} +13x_{9}+15x_{10} {}\leq 45.
\end{displaymath}
We then introduce the subsets 
$$N_1=\{1,2,3,4\},\ N_2=\{5\},\ \ldots,\ N_{7}=\{10\}.$$
We define 
\begin{align}
V_1 = \conv \Bigl\{\, (\ve x,\ve y^{N_1}) \in \Z_+^{4}\times \B^{4}:{}& 
3x_1+3x_2+3x_3+3x_4 =  \notag\\
&3y^{N_1}_1+6 y^{N_1}_2 + 9y^{N_1}_3 + 12y^{N_1}_4\notag\\ 
& y^{N_1}_1+y^{N_1}_2 + y^{N_1}_3 + y^{N_1}_4 \leq 1 \notag\\ 
& 0\leq x_i\leq 1, \;i=1,\ldots, 4 \qquad \qquad \Bigr\}. 
\end{align} 
Since $V_2,\ldots,V_7$ consist of single points each, these polyhedra are trivial. No additional $y$-variables are needed.
Then, $Q$ becomes 
\begin{align} 
Q = \conv\Bigl\{\,(\ve y^{N_1},x_5,\dots,x_{10})\in\B^{10}:{}&3y^{N_1}_1+6 y^{N_1}_2 + 9y^{N_1}_3 + 12y^{N_1}_4+4x_{5} \notag\\
&\quad{}+ 7x_{6} +8x_{7} +9x_{8} +13x_{9}+15x_{10}\leq 45\notag\\ 
& y^{N_1}_1+y^{N_1}_2 + y^{N_1}_3 + y^{N_1}_4 \leq 1,\notag\\
& x_i \in \{0,1\} \text{ for } i=5,\ldots,10 \quad 
\Bigr\}. 
\end{align} 
In this example, there are several other ways to define an extended formulation based on introducing new variables for the values that $x_1+x_2+x_3+x_4$ can attain.
One could introduce one particular integer variable $z$ that represents the value of $x_1+x_2+x_3+x_4$. Alternatively, one could introduce a binary expansion for the  values of  $x_1+x_2+x_3+x_4$, i.e., one introduces binary variables $z_1,z_2,z_3$ and requires that 
$ x_1+x_2+x_3+x_4 = z_1+2z_2+4z_3.$ For each of these models we compute the
facet description of the corresponding convex hull, as indicated in Table~\ref{tab:example-facets}. 
\begin{table}[b]\caption{Sizes of facet descriptions of various reformulations}
  \label{tab:example-facets}
  \begin{center}
  \begin{tabular}{lll}
    \toprule
    Formulation & Equations & \# Facets \\
    \midrule
    \mbox{original} &   & {}328  \\[1ex]
    \mbox{integer expansion} & $x_1+x_2+x_3+x_4 = z$  & {}328 \\[1ex]
    \mbox{binary expansion} &  $x_1+x_2+x_3+x_4 = z_1+2z_2+4z_3$  & {}217 \\[1ex]
    \mbox{value disjunction} & $x_1+x_2+x_3+x_4 = z_1+2z_2+3z_3+4z_4$  & {}77  \\
    & $z_1+z_2+z_3+z_4 \leq 1$\\
    \bottomrule
  \end{tabular}
\end{center}
\end{table}

In the original formulation there are $328$ facets needed to describe the polyhedron. If we introduce one additional integer variable that encodes the value of the  constraint $x_1+x_2+x_3+x_4$, then the same number of inequalities suffice to describe the corresponding convex hull of solutions. This is geometrically clear because every inequality of the original formulation is in bijection with an inequality in the lifted space. 
However, introducing three new binary variables $z_1,z_2,z_3$ and encoding the values of the partial constraint $x_1+x_2+x_3+x_4$ through the additional three variables $2^0z_1$, $2^1z_2$, $2^2z_3$, we obtain a polyhedron in the  $13$-dimensional space that requires $217$ facets for a complete description.
The value disjunction based on $x_1+x_2+x_3+x_4$ requires to introduce four new binary variables that are linked to the original variables by the two constraints
$$z_1+z_2+z_3+z_4 \leq 1,\quad x_1+x_2+x_3+x_4 = z_1+2z_2+3z_3+4z_4.$$
This new formulation in the $14$-dimensional space requires only $77$ facets
for a complete description.
\end{example} 

\begin{theorem}[Structure Theorem for Value Disjunction] 
\begin{equation} 
  \begin{aligned}
    P = \Bigl\{\,(\ve x,\ve w)\in \R^n\times \R^d:{}& \text{there exists }\ve y\in
    [0,1]^{n_1+\cdots+n_K} \text{ with }\\& (\ve y,\ve w) \in Q \text{ and }
    (\ve x^{N_i},\ve y^{N_i}) \in V_i \text{\emph{ for all }} i \,\Bigr\}. 
  \end{aligned} 
\label{claim to prove} 
\end{equation} 
\label{structure theorem for value disjunction} 
\end{theorem} 
\begin{proof}
The inclusion $\subseteq$ is trivial.
We shall prove the inclusion $\supseteq$. Let us consider $(\ve x,\ve w)$ from the set in 
the right-hand side of \eqref{claim to prove}. We try to prove that $(\ve
x,\ve w)\in P.$ 
For such an $(\ve x,\ve w)$, we know that there exists $\ve y$ such that $(\ve
x^{N_i},\ve y^{N_i}) \in V_i.$ Therefore there exist convex multipliers
$\lambda^{N_i,l} \geq 0$ with $\sum_{l=1}^{L_i} \lambda^{N_i,l} = 1$ such that 
\begin{equation} 
(\ve x^{N_i},\ve y^{N_i}) = \sum_{l=1}^{L_i} \lambda^{N_i,l} (\ve{\bar x}^{N_i,l},\ve{\bar 
y}^{N_i,l}), 
\label{combination Vi} 
\end{equation} 
where $(\ve{\bar x}^{N_i,l},\ve{\bar y}^{N_i,l})$ is an integral element of
$V_i$ and
$A(\ve{\bar y}^{N_i,l}) = A(\ve{\bar x}^{N_i,l})$. 
In particular the $y$-part is made of exactly one 1-entry. 
Therefore 
\begin{equation} 
y^{N_i}_t = \sum_{l\in T(N_i,t)} \lambda^{N_i,l} 
\label{y decomposition} 
\end{equation} 
with the sets $T(N_i,t)$, $t=1,\ldots, n_i$, being a packing of $\{1,\ldots, 
L_i\}$, namely for all $i$ we have
\begin{equation}
\{1,\ldots, L_i\} = T(N_i,1) \dotcup \ldots \dotcup T(N_i,n_i),
\label{packing using T}
\end{equation} 
where $C=A\dotcup B$ means $C=A\cup B$ and $A\cap B = \emptyset.$
The insight of \eqref{y decomposition} is shown in Figure \ref{fig y}. 
\begin{figure}[htpb] 
\begin{equation*} 
\left( 
\begin{array}{c} 
y^{N_1}_1 \\ \vdots \\ y^{N_1}_{n_1} \\ \hline \vdots \\ \hline 
y^{N_K}_1 \\ \vdots \\ y^{N_K}_{n_K} 
\end{array} 
\right) = \left( 
\begin{array}{c} 
\lambda^{N_1,\cdot} + \cdots + \lambda^{N_1,\cdot}\\ 
\vdots \\ 
\lambda^{N_1,\cdot} + \cdots \\ 
\hline 
\vdots \\ 
\hline 
\lambda^{N_K,\cdot} + \cdots \\ 
\vdots \\ 
\lambda^{N_K, \cdot} + \cdots 
\end{array} 
\right) 
\end{equation*} 
\caption{Each $y$ is equal to the sum of zero, one or more  $\lambda$ 
from the convex combination.} 
\label{fig y} 
\end{figure} 
\medskip 
 
Up to now we have used the fact that $(\ve x^{N_i},\ve y^{N_i}) \in V_i.$ 
We also have a second condition stating that $(\ve y,\ve w) \in Q.$ Therefore there exist 
convex multipliers $\sigma_r \geq 0$ with $\sum_{r=1}^R \sigma_r = 1$ such that 
\begin{equation} 
\ve y = \sum_{r=1}^R \sigma_r \ve{\hat y}^r \quad\text{and}\quad \ve w=\sum_{r=1}^R \sigma_r \ve{\hat w}^r, 
\label{combination Q} 
\end{equation} 
where 
\begin{equation*} 
\ve{\hat y}^r = {\bigl(\ve{\hat y}^{N_1,r}, \ldots, \ve{\hat y}^{N_K,r}\bigr)}, 
\end{equation*} 
and where $\ve{\hat y}^{N_i,r}$ is a unit vector. Furthermore 
\begin{equation*} 
\sum_{i=1}^K A(\ve{\hat y}^{N_i,r})+\sum_{j=1}^d G_j \ve{\hat w}^r_j \leq \ve b.
\end{equation*} 

We are now able to express $(\ve x,\ve w)$ as a convex combination of feasible 
solutions of $A\ve x+G\ve w\leq \ve b$, using the convex combinations \eqref{combination Q} 
and \eqref{combination Vi}. 
To do this, we first remark that, similarly to \eqref{y decomposition}, we can 
express $\ve y$ in terms of $\sigma_r$ only, namely 
\begin{equation} 
y^{N_i}_t = \sum_{s\in S(N_i,t)} \sigma_s, 
\label{y sigma decomposition} 
\end{equation} 
with the sets $S(N_i,t)$, $t=1,\ldots, n_i$ being a packing of $\{1,\ldots, 
R\},$ namely
\begin{equation}
\{1,\ldots, R\} = S(N_i,1) \dotcup \ldots \dotcup S(N_i,n_i),
\label{packing using S}
\end{equation}
for all $i$.
By using \eqref{y decomposition}, we therefore conclude that 
\begin{equation} 
\sum_{s\in S(N_i,t)} \sigma_s = \sum_{l\in T(N_i,t)} \lambda^{N_i,l} 
\label{linking lambda sigma} 
\end{equation} 
By using the similarity of decompositions \eqref{y sigma decomposition} and 
\eqref{y decomposition}, we can construct the desired convex 
combination as follows.\medskip 
 
Let us fix $r$, i.e., we consider each pair $(\sigma_r,\ve{\hat y}^r)$ separately.  We know that $\ve{\hat y}^r$ is 
divided into $K$ blocks with a unit vector in each block. In the block $N_i$, 
we refer to the index of the non-zero component of $\ve{\hat y}^r$ as $c(\ve{\hat y}^{N_i,r}).$ 
Using 
\eqref{y decomposition}, we can associate to $c(\ve{\hat y}^{N_i,r})$ a set $T(N_i, 
c(\ve{\hat y}^{N_i,r}))$ of indices~$l$, which correspond to multipliers
$\lambda^{N_i,l}$ and vectors $\ve{\bar x}^{N_i, l}$ of the convex combination~\eqref{combination Vi}. 
For every possible choice of indices 
\begin{displaymath}
  l^r_1\in T(N_1, c(\ve{\hat y}^{N_1,r})),\quad \ldots,\quad l^r_K\in T(N_K,c(\ve{\hat 
y}^{N_K,r})),
\end{displaymath}
we consider the point
\begin{equation*} 
\ve x(l^r_1,\ldots, l^r_K)={\bigl(\ve{\bar x}^{N_1,l^r_1},  
\cdots, \ve{\bar x}^{N_K,l^r_K}\bigr)} 
\end{equation*} 
with a corresponding coefficient
\begin{equation} 
  \nu(l^r_1,\ldots, l^r_K)= 
\sigma_r \frac{\lambda^{N_1,l^r_1}}{\displaystyle\sum_{l\in T(N_1,c(\ve{\hat y}^{N_1, 
r}))} 
\lambda^{N^1,l}} \cdots \frac{\lambda^{N_K,l^r_K}}{\displaystyle\sum_{l\in T(N_K 
,c(\ve{\hat 
y}^{N_K,r}))} 
\lambda^{N^K,l}}. 
\label{definition nu} 
\end{equation} 
First we can see that for all $l^r_1, \ldots, l^r_K$, the vector 
$(\ve x(l^r_1,\ldots, l^r_K),\ve{\hat w}^r)$ satisfies $A\, \ve x(l^r_1,\ldots, l^r_K)+G\ve{\hat w}^r \leq \ve b$. 
Indeed,
\begin{align*} 
A\,\ve x(l^r_1,\ldots, l^r_K) +G\ve{\hat w}^r& = A(\ve{\bar x}^{N_1,l^r_1}) + \cdots + A(\ve{\bar 
x}^{N_K,l^r_K})+G\ve{\hat w}^r\\ 
& = A(\ve{\hat y}^{N_1,l^r_1}) + \cdots + A(\ve{\hat y}^{N_K,l^r_K})+G\ve{\hat w}^r\\ 
& = A(\ve{\hat y}^r) +G\ve{\hat w}^r\\ 
&\leq \ve b, 
\end{align*} 
since $(\ve{\hat y}^r,\ve{\hat w}^r)$ is a mixed-0/1 solution of~$Q$. It
now suffices to prove that $\ve x$ 
is the convex combination of all the $\ve x(l^r_1,\ldots, l^r_K)$ using the 
corresponding 
coefficients $\nu(l^r_1,\ldots, l^r_K).$ 
Let us fix $N_i$ and an index $j\in N_i$. We have 
\begin{align} 
x^{N_i}_j &=\sum_{r=1}^R \sum_{l_1^r\in T(N_1,c(\ve{\hat y}^{N_1,r}))} \cdots 
\sum_{l_K^r \in T(N_K,c(\ve{\hat y}^{N_K,r}))} \nu(l_1^r,\ldots, l_K^r) 
x^{N_i}_j(l_1^r, \ldots, l_K^r)\notag\\ 
&=\sum_{r=1}^R \sum_{l_1^r\in T(N_1,c(\ve{\hat y}^{N_1,r}))} \cdots 
\sum_{l_K^r \in T(N_K,c(\ve{\hat y}^{N_K,r}))} \nu(l_1^r,\ldots, l_K^r) \bar 
x^{N_i,l_i^r}_j\notag\\ 
&=\sum_{r=1}^R \sum_{l_i^r\in T(N_i, c(\ve{\hat y}^{N_i,r}))}\sigma_r 
\frac{\lambda^{N_i,l_i^r}}{\displaystyle\sum_{l\in T(N_i,c(\ve{\hat y}^{N_i,r}))} 
\lambda^{N_i,l}} \bar x^{N_i,l_i^r}_j,\label{last line} 
\end{align} 
the last identity being obtained using \eqref{definition nu}. 
For a fixed $i$, we have, using \eqref{packing using S},
\begin{equation*}
\{1,\ldots,R\} = S(N_i,1) \dotcup \ldots \dotcup S(N_i,n_i).
\end{equation*}
Therefore we can rewrite \eqref{last line} using indices running over the
different $S(N_i,k).$ Remark also that when we fix $r\in S(N_i,k)$, we have
$c(\hat{\ve y^{N_i,r}}) = k.$ We hence have
\begin{align}
x_j^{N_i} &= \sum_{k=1}^{n_i} \sum_{p\in S(N_i,k)} \sum_{l\in T(N_i,k)} 
\sigma_p \frac{\lambda^{N_i,l}}{\displaystyle\sum_{q\in T(N_i,k)} \lambda^{N_i,q}} \bar
x_j^{N_i,l}\notag\\
&= \sum_{k=1}^{n_i} \sum_{l\in T(N_i,k)} \frac{\displaystyle\sum_{p\in S(N_i,k)}
\sigma_p}{\displaystyle\sum_{q\in T(N_i,k)} \lambda^{N_i,q}} \lambda^{N_i,l} \bar
x_j^{N_i,l}\notag\\
&= \sum_{k=1}^{n_i} \sum_{l\in T(N_i,k)} \lambda^{N_i,l} \bar x_j^{N_i,l},
\label{sum over k,l}
\end{align}
where \eqref{sum over k,l} is obtained using \eqref{linking lambda sigma}.
We can use \eqref{packing using T} namely
\begin{equation*}
T(N_i,1) \dotcup \ldots \dotcup T(N_i,n_i) = \{1,\ldots, L_i\}.
\end{equation*}
In particular it allows us to sum over $\{1,\ldots, L_i\}$ in \eqref{sum over
k,l} instead of the summation over $k$ and $l$. We therefore finally have
\begin{equation*}
x_j^{N_i} = \sum_{l=1}^{L_i} \lambda^{N_i,l} \bar x_j^{N_i,l},
\end{equation*}
which is the desired result using \eqref{combination Vi}.
Finally, the sum of the $\nu$ coefficients is equal to 1 due to their construction and 
the fact that $\sum_{r=1}^R \sigma_r =1.$
\end{proof}

\begin{example}\label{ex:fourvariable}
Consider the set 
\begin{equation*} 
\mathcal F=\{x\in\{0,1,2\}^4: x_1+x_2+2x_3+3x_4 \leq 7\}. 
\end{equation*} 
The complete facet description of $\conv \mathcal F$ is given by the 14 
inequalities $\ve c^\top \ve x \leq \gamma$ shown in Table~\ref{initial description}. 
\begin{table}[htbp] 
\caption{The complete description of Example~\ref{ex:fourvariable} in the original space} 
\label{initial description} 
\begin{center} 
  \small
  $\begin{array}{*4{r}@{\;}lc@{\qquad}*4{r}@{\;}l}
    \toprule
    c_1 & c_2 & c_3 & c_4 & \multicolumn{1}{r}{\gamma}
    && c_1 & c_2 & c_3 & c_4 & \multicolumn{1}{r}{\gamma}\\
    \cmidrule{1-5}\cmidrule{7-11}
    -1 & 0 & 0 & 0 &{}\leq 0 &&  1 & 0 & 0 & 1 &{}\leq  3\\  
    0 & -1 & 0 & 0 &{}\leq 0 &&  0 & 1 & 0 & 1 &{}\leq  3\\  
    0 & 0 & -1 & 0 &{}\leq 0 &&  0 & 0 & 1 & 2 &{}\leq  4\\ 
    0 & 0 & 0 & -1 &{}\leq 0 &&  1 & 1 & 1 & 1 &{}\leq  5\\ 
    1 & 0 & 0 & 0 &{}\leq 2  &&  0 & 1 & 2 & 2 &{}\leq  6\\ 
    0 & 1 & 0 & 0 &{}\leq 2  &&  1 & 0 & 2 & 2 &{}\leq  6\\ 
    0 & 0 & 1 & 0 &{}\leq 2  &&  1 & 1 & 2 & 3 &{}\leq  7\\ 
    \bottomrule
  \end{array}$
\end{center}
\end{table} 
 
We now construct a value disjunction of the set $\mathcal F$. 
To do this, we consider three blocks $N_1=\{1,2\}$, $N_2=\{3\}$, $N_4=\{4\}$. 
In block $N_1$ we consider the linear form $x_1+x_2$, which can take the
values $0,1,\dots,4$ because $x_1$ and $x_2$ have an upper bound of~$2$.
We introduce thus four 
variables $y_1,y_2,y_3,y_4$ corresponding to the four nonzero values. The blocks $N_2$ 
and $N_3$ are trivial, so we do not need to introduce new variables in those 
cases. A valid 
formulation for $\mathcal F$ is thus 
\begin{align*} 
\mathcal F = \Proj_{\ve x}\bigl\{\,(\ve x,\ve y)\in \{0,1,2\}^4\times\B^4:{} 
& y_1+2y_2+3y_3+4y_4+2x_3+3x_4 
\leq 7\\ 
& x_1+x_2 = y_1+2y_2+3y_3+4y_4\\ 
& y_1+y_2+y_3+y_4 \leq 1\,\bigr\}. 
\end{align*} 
Theorem \ref{structure theorem for value disjunction} now asserts that 
we obtain the complete description of the extended formulation of~$\mathcal F$
by combining the complete descriptions of the polyhedra
\begin{align*} 
V_1 = \conv\{\,(x_1,x_2,\ve y)\in \{0,1,2\}^2\times\B^4: {}
& x_1+x_2 = y_1+2y_2+3y_3+4y_4\\ 
& y_1+y_2+y_3+y_4 \leq 1\qquad\qquad \}, 
\end{align*} 
and 
\begin{align*} 
Q = \conv\{\,(x_3,x_4,\ve y) \in \{0,1,2\}^2\times \B^4: {}
& 2x_3+3x_4 + y_1+2y_2+3y_3+4y_4 \leq 7 \\ 
&y_1+ y_2+ y_3 + y_4 \leq 1 \qquad \qquad\}. 
\end{align*} 
We obtain the facet description given 
by the inequalities $\ve c^\top\ve x + \ve d^\top\ve y \leq \gamma$ shown 
in Table~\ref{description in extended space}. 
For each non-trivial inequality, we also mention whether it comes from $V_1$ 
or from $Q$. 
\begin{table}[htbp] 
\caption{The complete description of Example~\ref{ex:fourvariable} in the extended space}
\label{description in extended space} 
\small
\begin{center}
  $\begin{array}{*8{r}@{\;}ll}
    \toprule
    c_1 & c_2 & c_3 & c_4 & d_1 & d_2 & d_3 & d_4 & \multicolumn{1}{r}{\gamma}
    & \text{ Origin }\\ 
    \midrule 
      &   & -1 &   &   &   &   &   &{}\leq 0 \\ 
      &   &   & -1 &   &   &   &   &{}\leq 0 \\ 
      &   &   &   & -1 &   &   &   &{}\leq 0 \\ 
      &   &   &   &   & -1 &   &   &{}\leq 0 \\ 
      &   &   &   &   &   & -1 &   &{}\leq 0 \\ 
      &   &   &   &   &   &   & -1 &{}\leq 0 \\ 
      & -1 &   &   &   &   & 1 & 2 &{}\leq 0 & V_1\\ 
      & -1 &   &   & 1 & 2 & 2 & 2 &{}\leq 0 & V_1\\ 
      &   &   &   & 1 & 1 & 1 & 1 &{}\leq 1 & Q, V_1\\ 
      &   & 1 &   &   &   &   & 1 &{}\leq 2 & Q\\ 
      &   &   & 1 &   & 1 & 1 & 1 &{}\leq 2 & Q\\ 
      &   & 1 & 1 & 1 & 1 & 1 & 2 &{}\leq 3 & Q\\ 
      &   & 1 & 2 &   & 1 & 2 & 2 &{}\leq 4 & Q\\ 
    1 & 1 &   &   & -1 & -2 & -3 & -4 & {}= 0 & V_1\\ 
    \bottomrule
  \end{array}$
\end{center}
\end{table} 
\end{example} 
 
In the example it turns out that the number of inequalities describing $\conv
\mathcal F$ is 
the same in the two representations. This, however, is not always true.  Moreover, an inherent advantage of the second 
formulation over the first formulation is that its structure is better known. 
In particular, it may occur that the same polyhedron $V_i$ appears in several 
different problems. In this case, the knowledge about the description of the polyhedron  $V_i$ can be used over and over again. 

The next section presents the case of a polyhedron that 
appears often in our experiments, namely the $V_i$ polyhedron where all the 
coefficients of the variables $x$ are the same. We show that we can compute 
a full description for this object.

\section{A special family of linking polyhedra} 
 
 
In this section we study the linking polyhedra $V_i$ for the case where the 
columns $A_j$ for $j\in N_i$ are identical and the variables $x_j$ are binary. 
In other words, we study the polytope 
\begin{align*} 
  \label{eq:cardinality-counting} 
  V_i=\conv \{(\ve x^{N_i},\ve y^{N_i}) \in \B^{|N_i|}\times \B^{n_i}:& 
  \sum_{j\in N_i} x_j = \sum_{k=1}^{n_i} k y^{N_i}_k\notag\\ 
  & \sum_{k=1}^{n_i} y^{N_i}_k \leq 1 \qquad \qquad \}. 
\end{align*} 
We are able to give a complete description of this polytope $V_i$. 
 
\begin{theorem} \label{th:cardinality-counting-facets}
  \begin{subequations} 
    \label{eq:cardinality-counting-facets} 
    $V_i$ is a polytope whose affine hull is given by the equation: 
    \begin{equation} 
      \sum_{j\in N_i} x_j = \sum_{k=1}^{n_i} k y^{N_i}_k \label{eq:cardinality-counting-subspace} 
    \end{equation} 
    The facets of $V_i$ are given by: 
    \begin{align} 
      \sum_{j\in T} x_j - \sum_{k=1}^{|T|} k y_k - \sum_{k=|T|+1}^{n_i} |T| y_k 
      &\leq 0 
      && \text{for $\emptyset\neq T\subset N_i$}\label{eq:cardinality-counting-facet-a}\\ 
      \sum_{k=1}^{n_i} y^{N_i}_k &\leq 1 \label{eq:cardinality-counting-facet-b}\\ 
      y^{N_i}_k &\geq 0 && \text{for $k=1,\dots,n_i$.} \label{eq:cardinality-counting-facet-c} 
    \end{align} 
  \end{subequations} 
\end{theorem} 
\begin{proof} 
  We first show that the inequalities~\eqref{eq:cardinality-counting-facets} 
  are valid for $V_i$.  To this end, let $(\ve x, \ve y) \in\B^{|N_i|}\times 
  \B^{n_i}$ be a vertex of $V_i$.  If $\ve y = \ve0$, then also $\ve x=\ve0$, 
  and inequality~\eqref{eq:cardinality-counting-facet-a} is trivially satisfied. 
  Otherwise, $\ve y = \ve e^k$ with $k = \sum_{j\in N_i} x_j = |\supp \ve x^{N_i}|$.  
  Let $\emptyset\neq T\subset N_i$ be 
  arbitrary.  If $k\leq |T|$, we have 
  \begin{displaymath} 
    \sum_{j\in T} x_j - \sum_{k=1}^{|T|} k y_k - \sum_{k=|T|+1}^{n_i} |T| y_k 
    = \sum_{j\in T} x_j - k \leq 0. 
  \end{displaymath} 
  On the other hand, if $k > |T|$, we have 
  \begin{displaymath} 
    \sum_{j\in T} x_j - \sum_{k=1}^{|T|} k y_k - \sum_{k=|T|+1}^{n_i} |T| y_k 
    = \sum_{j\in T|} x_j - |T| \leq 0. 
  \end{displaymath} 
  Hence, \eqref{eq:cardinality-counting-facet-a} is satisfied.  
  The remaining inequalities are trivially valid for $V_i$. 
 
  For the ease of  
  notation we let $N=N_i$, $n=|N|$ and substitute the variables $y^{N_i}_k$ by  
  simply $y_k$. Let $\ve c^\top\ve x+\ve d^\top\ve y \leq \gamma$ be a facet-defining inequality of  
  $V_i$ and set 
  $$F=\{\,(\ve x,\ve y) \in V_i : \ve c^\top\ve x+\ve d^\top\ve y = \gamma\,\}.$$ 
  We will show that $\ve c^\top\ve x+\ve d^\top\ve y \leq \gamma$ corresponds to one of the inequalities  
  in~\eqref{eq:cardinality-counting-facets} up to multiplication by a scalar. We assume that the variables in $N$ are  
  reordered such that  
  $c_1\geq c_2 \geq \ldots \geq c_n$. Since $V_i$ is not full dimensional, we first  
  transform $\ve c^\top\ve x+\ve d^\top\ve y \leq \gamma$ into a standard form. This can be achieved by  
  adding multiples of the equation~\eqref{eq:cardinality-counting-subspace} 
  to $\ve c^\top\ve x+\ve d^\top\ve y \leq \gamma$. 
  More precisely, we first proceed with the following two steps. 
  \begin{enumerate}[(1)] 
  \item While there exists an index $i \in N$ such that $c_i<0$, add $-c_i$  
  times Equation~\eqref{eq:cardinality-counting-subspace} to the inequality 
  $\ve c^\top\ve x+\ve d^\top\ve y \leq \gamma$. 
  Let us again denote by $\ve c^\top\ve x+\ve d^\top\ve y \leq \gamma$ the resulting inequality.  
  Notice that  after terminating with Step 1, we have that  $c_i \geq 0$ for all  
  $i \in N$ and $c_n=0$. 
  \item If $c_i >0$ for all $i \in N$ and there exist $i,j \in N$ such that  
    $c_i \neq c_j$, then $c_1 > c_n >0$ due to our reordering. In this case we  
  subtract $c_n$ times  Equation~\eqref{eq:cardinality-counting-subspace} from the inequality  
  $\ve c^\top\ve x+\ve d^\top\ve y \leq \gamma$. Notice that also after Step (2) has 
  been performed we   
  have that $c_n=0$ and $c_i \geq 0$ for all $i \in N$.  
  \end{enumerate} 
  The preprocessing steps (1) and (2) guarantee that $c_i \geq 0$ for all $i \in  
  N$.  Now let $s\in\{0,\dots,n\}$ be an index such that  
   $$c_1\geq c_2 \geq \ldots\geq c_s> 0 = c_{s+1} = \ldots = c_n.$$ 
  We define $T=\{\,i \in N : c_i >0\,\} = \{1,\dots,s\}$. 
  We consider the following cases. 
 
  \begin{enumerate}[{Case} 1.] 
  \item If $T=\emptyset$, i.e., $c_1=\dots=c_n=0$, it follows that $\ve c^\top\ve 
    x+\ve d^\top\ve y \leq \gamma$ is a multiple of the inequality $\sum_{k=1}^n 
    y_k \leq 1$ or of the non-negativity constraints $y_k \geq 0$.  
     
    Indeed, because $(\ve0,\ve0)$ is feasible, we have $\gamma\geq0$.  Since 
    $F$ is a facet, there must be $2n-1$ affinely independent feasible points 
    on it. If $\gamma=0$, we have $(\ve 0,\ve0)\in F$; therefore, for all but 
    one $k=1,\dots,n$, a point $(\ve x, \ve e^k)$ must be contained in $F$. 
    This means that $d_k=\gamma=0$ for all but one $k=1,\dots,n$.  For the 
    remaining one $\tilde k\in\{1,\dots,n\}$ we have $d_{\tilde k}\leq \gamma=0$, 
    so $\ve c^\top \ve x + \ve d^\top \ve y \leq \gamma$ is a scalar multiple of the 
    non-negativity constraint $y_{\tilde k} \geq0$. 
 
    On the other hand, if $\gamma>0$, then $(\ve0,\ve 0)\notin F$, so we have 
    $F\subseteq\{\,(\ve x,\ve y) \in  V_i : \sum_{k=1}^n y_k =1\,\}$, since 
    $(\ve0,\ve0)$ is the only feasible integer point with $\ve y=\ve0$. 
    Because $F$ is a facet, we have $F = \{\,(\ve x,\ve y) \in  V_i : 
    \sum_{k=1}^n y_k =1\,\}$, which corresponds 
    to~\eqref{eq:cardinality-counting-facet-b}.  
 
  \item If $T=N$, we conclude from our previous analysis that $c_i = c_j \neq 
    0$ for all $i,j \in N$. It follows that $\ve c^\top\ve x+\ve d^\top\ve y \leq 
    \gamma$ is implied by Equation~\eqref{eq:cardinality-counting-subspace}, a 
    contradiction that $F$ defines a facet of $V_i$. 
 
  \item Therefore, we may assume that $\emptyset\neq T\subset N$, $T \neq N$. 
    Again, since $(\ve0,\ve0)$ is feasible, we have that $\gamma \geq 
    0$. If $\gamma >0$, then $F \subseteq \{\,(\ve x,\ve y) \in  V_i : 
    \sum_{k=1}^n y_k =1\}$. 
    Hence, we can assume that $\gamma = 0$.  
 
    We next define indices $1\leq i_1<i_2<\ldots<i_r\leq s$ as follows:    
    $$c_1=\ldots =c_{i_1} > c_{i_1+1} = \ldots = c_{i_2}> \ldots > c_{i_r+1} =  
    \ldots = c_{s}.$$ 
    By testing the inequality $\ve c^\top\ve x+ \ve d^\top\ve y\leq 0$ 
    with the feasible points $(\ve e^1, \ve e^1)$, $(\ve e^1+\ve e^2, \ve e^2)$, 
    $(\ve e^1+\ve e^2+\ve e^3, \ve e^3)$, \dots, we conclude that 
    \begin{align*} 
      -d_1 & \geq  c_1\\ 
      -d_2 & \geq  c_1+c_2\\ 
      &\vdots\\ 
      -d_{i_1} & \geq  c_1+c_2+\ldots+c_{i_1}\\ 
      -d_{i_1+1} & \geq  \textstyle\sum_{j=1}^{i_1} c_j + c_{i_1+1}\\ 
      -d_{i_1+2} & \geq  \textstyle\sum_{j=1}^{i_1} c_j + c_{i_1+1}+c_{i_1+2}\\ 
      &\vdots\\ 
      -d_{i_2} & \geq \textstyle\sum_{j=1}^{i_1} c_j +c_{i_1+1}+c_{i_1+2}+\ldots+c_{i_2}\\ 
      &\vdots\\ 
      -d_{i_r+1} & \geq  \textstyle\sum_{j=1}^{i_r} c_j + c_{i_r+1}\\ 
      -d_{i_r+2} & \geq  \textstyle\sum_{j=1}^{i_r} c_j + c_{i_r+1}+c_{i_r+2}\\ 
      &\vdots\\ 
      -d_{s} & \geq \textstyle\sum_{j=1}^{i_r} c_j +c_{i_r+1}+c_{i_r+2}+\ldots+c_{s} 
    \end{align*} 
    Therefore, the inequality $\ve c^\top\ve x+\ve d^\top\ve y \leq \gamma=0$ is 
    dominated by the  
    following conic combination of the 
    inequalities~\eqref{eq:cardinality-counting-facet-a}:  
    $$\begin{array}{r@{}c@{}l} 
      c_{i_r} & {}\times{} & \biggl( \sum_{i=1}^{s_{}} x_i - \sum_{k=1}^s ky_k -  
      \sum_{k=s_{}+1}^n sy_k \leq 0 \biggr)\\ 
      {}+ (c_{i_r}-c_{i_r-1}) & \times & \biggl( \sum_{i=1}^{i_r} x_i - \sum_{k=1}^{i_r}  
      ky_k - \sum_{k=i_r+1}^n i_r y_k \leq 0 \biggr) \\ 
      &\vdots\\ 
      {}+ (c_{i_1}-c_{i_2}) & \times & \biggl( \sum_{i=1}^{i_1} x_i - \sum_{k=1}^{i_1}  
      ky_k - \sum_{k=i_1+1}^n i_1 y_k \leq 0 \biggr). \\ 
    \end{array}$$ 
  \end{enumerate} 
  This completes the proof. 
\end{proof} 
 
\begin{theorem} 
  The separation problem over the linking polyhedron $V_i$ in the case of
  identical coefficients can be solved in polynomial time. 
\end{theorem} 
\begin{proof} 
  Let $(\ve x^*,\ve y^*)$ be a point satisfying the polynomially many 
  constraints~(\ref{eq:cardinality-counting-subspace}, 
  \ref{eq:cardinality-counting-facet-b}, 
  \ref{eq:cardinality-counting-facet-c}).  We show that, in polynomial time, we 
  can decide whether $(\ve x^*,\ve y^*)$ satisfies the exponentially many 
  inequalities~\eqref{eq:cardinality-counting-facet-a}; if it does not, we 
  can construct a maximally violated inequality. 
 
  It is clear that among the inequalities 
  \eqref{eq:cardinality-counting-facet-a} with equal cardinality~$|T| = s$, a 
  most violated inequality is the one where $T$ is the index set of the $s$ 
  largest components $x_j^*$.  Therefore it suffices to sort the variables 
  $x_1^*,\dots,x^*_{|N_i|}$ such that  
  \begin{displaymath} 
    x_1^*\geq x_2^*\geq\dots\geq x^*_s>0=x_{s+1}=\dots=x^*_{|N_i|}. 
  \end{displaymath} 
  Then we can simply evaluate the violation of 
  inequality~\eqref{eq:cardinality-counting-facet-a} for the sets $\{1\}$, 
  $\{1,2\}$, $\{1,2,3\}$, \dots, $\{1,\dots,s\}$ and pick the set which yields 
  the maximal violation.  
\end{proof} 
 
\section{An application: The knapsack with three distinct coefficients} 
 
In this section, we show that the value disjunction procedure is a tool to 
compute complete descriptions in an extended space. As an example  
we consider the 0/1 knapsack problem with three distinct 
coefficients: 
\begin{equation} 
  \label{eq:knapsack-3coeff} 
  \sum_{j\in N_1} \mu x_j + \sum_{j\in N_2} \lambda x_j + \sum_{j\in N_3} 
  \sigma x_j \leq \beta, 
\end{equation} 
where $N_1$, $N_2$, $N_3$ are pairwise disjoint index sets.  The convex hull 
of the feasible solutions can have exponentially many vertices and facets. 
Moreover, the complete 
facet description for~\eqref{eq:knapsack-3coeff} is not known in general.   
In \cite{Weismantel96}, the case of the knapsack with two different 
coefficients was solved. 
By 
applying the structure theorem for value disjunctions 
(Theorem \ref{structure theorem for value disjunction}), we are able to give a complete 
description for an  
extended formulation of~\eqref{eq:knapsack-3coeff} using only polynomially 
many variables.  
 
We consider the extended formulation of~\eqref{eq:knapsack-3coeff}, 
\begin{subequations}  
  \begin{align*}  
    \sum_{j\in N_1} \mu x_j + \sum_{j\in N_2} \lambda x_j + \sum_{j\in 
      N_3}  
    \sigma x_j &\leq \beta\\ 
    \sum_{j\in N_i} x_j &= \sum_{k=1}^{|N_i|} k y^i_k  & & \text{for $i=1,2,3$} \\ 
    \sum_{k=1}^{|N_i|} y^i_k &\leq 1 & & \text{for $i=1,2,3$} \\ 
    x&\in\B^{|N_1| + |N_2| + |N_3|}\\ 
    y^i &\in\B^{|N_i|} & & \text{for $i=1,2,3$}.  
  \end{align*}  
\end{subequations}  
Theorem~\ref{structure theorem for value disjunction} provides us the
framework to describe the convex hull of such an extended formulation.  It is
given by the intersection of the linking polyhedron and the aggregated
polyhedron.  The linking polyhedron was studied in the last section.
Theorem~\ref{th:cardinality-counting-facets} gives a complete facet
description of it. Concerning the aggregated polyhedron, we will make use of
a vertex description.  It is the convex hull of the set described by  
\begin{subequations} 
  \begin{align*} 
    \mu \sum_{k=1}^{|N_1|} k y^{N_1,k} + \lambda \sum_{k=1}^{|N_2|} k y^{N_2,k}   
    + \sigma \sum_{k=1}^{|N_3|} k y^{N_3,k} &\leq \beta \\ 
    \sum_{k=1}^{|N_i|} y^{N_i,k} &\leq 1 & & \text{for $i=1,2,3$} \\ 
    \ve y^{N_i} &\in\B^{|N_i|} & & \text{for $i=1,2,3$}.  
  \end{align*} 
\end{subequations} 
Clearly there are at most $(1+|N_1|)\cdot (1+|N_2|)\cdot (1+|N_3|)$ vertices. 
We denote them by $\ve v^1,\dots,\ve v^p\in \B^{|N_1|+|N_2|+|N_3|}$. 
 
\begin{theorem} 
  The complete facet description of~\eqref{eq:knapsack-3coeff} 
  in an extended space is given by: 
  \begin{subequations} 
    \begin{align*} 
      &\ve y  = \sum_{j=1}^p \ve v^j z_j \\ 
      & \sum_{j=1}^p z_j = 1 \\ 
      &z_j \geq 0 && \text{for $j=1,\dots,p$} \\ 
      &\sum_{j\in N_i} x^{N_i}_j = \sum_{k=1}^{n_i} k y^{N_i,k} 
      && \text{for $i=1,2,3$} 
      \\ 
      &\sum_{j\in T} x^{N_i}_j \geq \sum_{\substack{k\in\{1,\dots,n_i\}:\\ |T| + k > n_i }} (|T| + k - n_i) y^{N_i,k} 
      && \text{for $i=1,2,3$ and $\emptyset\neq T\subset 
        N_i$} 
      \\ 
      &\ve x\in \R^{|N_1|+|N_2|+|N_3|}\\ 
      &\ve y\in \R^{|N_1|+|N_2|+|N_3|}\\ 
      &\ve z\in\R^p. 
    \end{align*}%
  \end{subequations} 
\end{theorem} 
\begin{proof} 
This follows from Theorem \ref{structure theorem for value disjunction}. 
\end{proof}

It is straightforward to extend our construction to binary integer programs
with a fixed number of different columns.  
 
\section{Branching on value disjunctions} 
 
So far we have presented the value disjunction technique as a theoretical
tool to define extended formulations which may yield more tractable
polyhedral descriptions.  Clearly it would be too much to expect general results on 
the existence or constructability of an intermediate representation
for an arbitrary integer program that is better than the original formulation.  
The more modest goal of this section is to provide evidence for the practical 
usefulness of the value disjunction technique, using a limited set of
computational experiments.  

We shall restrict ourselves to experiments where we perform branching on the
new binary variables of the extended formulation.  We first need to discuss
the situations for which we propose to make use of our new technique, so as to
complement the existing branch-and-cut techniques.

\paragraph{On the simplification effect of branching.}

Today mixed integer linear programs are solved using branch-and-cut algorithms, i.e., such an algorithm consists of two phases, 
the cutting phase with the objective to tighten a 
current formulation and a branching phase. 
However as of today there are essentially no 
mathematical arguments available that help to decide 
when it is more efficient to branch or to cut.
This question is fundamental since computational experiments 
clearly reveal that neither a pure branch-and-bound algorithm 
nor a pure cutting plane algorithm can solve the instances 
that the combination of the two can manage to solve. 
One partial answer to this question is given by the fact  that branching does
not only generate subproblems with less variables, but, more importantly, the
polyhedral description of each of the two subproblems is significantly easier than the original facet description. We illustrate this point through an example.

\begin{example}
  \label{ex:branching}
  We consider the feasible region
  \begin{align*}
    7x_1+5x_2-x_3-x_4-2x_5-3x_6-4x_7-6x_8 \leq 1 \\
    x_i\in\{0,1\}.
  \end{align*}
  The non-trivial facets of the convex hull are shown in
  Table~\ref{table:ex-full}.  If we consider the four subproblems where the
  variables $x_7$ and $x_8$ are fixed to the possible values, we obtain much
  simpler facet descriptions; see Table~\ref{table:ex-branches}.
\end{example}
\begin{table}[tbp]
    \caption{Full description of Example~\ref{ex:branching}}  \label{table:ex-full} 
    \begin{center}
      \small
      $\begin{array}{*8{r@{\ }}@{\;}lc@{\qquad}*8{r@{\ }}@{\;}l}
  \toprule
  c_1 & c_2 & c_3 & c_4 & c_5 & c_6 & c_7 & c_8 & \multicolumn{1}{c}{\gamma}
  && c_1 & c_2 & c_3 & c_4 & c_5 & c_6 & c_7 & c_8 & \multicolumn{1}{c}{\gamma}\\
  \cmidrule{1-9}\cmidrule{11-19}
   &{}    &{}   &{}   &-1&{}    &-1&-1        &{}\leq 0&&     3&2&-1&-1&-1&-1&-1&-2        &{}\leq 1\\     
   &{}    &{}   &{}   &{}    &-1&-1&-1        &{}\leq 0&&      2&2&-1&-1&-2&{}    &-2&-1        &{}\leq 1\\  
  &{}    &-1&-1&{}    &{}    &-1&-1        &{}\leq 0&&         3&2&-1&-1&-2&{}    &-2&-2        &{}\leq 1\\  
  &{}    &-1&{}   &-1&-1&{}    &-1        &{}\leq 0&&         3&1&-1&-1&-2&-2&{}    &-2        &{}\leq 1\\  
  &{}    &{}   &-1&-1&-1&{}    &-1        &{}\leq 0&&         3&3&-1&-1&-2&-1&-2&-2        &{}\leq 1\\    
  &1&-1&{}   &{}    &-1&-1&-1        &{}\leq 0&&           3&2&-1&{}   &-1&-1&-1&-3        &{}\leq 1\\  
  &1&{}   &-1&{}    &-1&-1&-1        &{}\leq 0&&           3&2&{}   &-1&-1&-1&-1&-3        &{}\leq 1\\  
  &1&{}   &{}   &-1&-1&-1&-1        &{}\leq 0& &          3&2&{}   &{}   &-1&-1&-2&-3        &{}\leq 1\\
  &1&-1&-1&-1&{}    &-1&-1        &{}\leq 0&   &          3&3&-1&{}   &-1&-2&-3&-2        &{}\leq 1\\  
 2&1&-1&{}   &-1&-1&-1&-2        &{}\leq 0&    &         3&3&{}   &-1&-1&-2&-3&-2        &{}\leq 1\\  
 2&1&{}   &-1&-1&-1&-1&-2        &{}\leq 0&    &         4&2&-1&{}   &-2&-1&-2&-3        &{}\leq 1\\  
 2&1&{}   &{}   &-1&-1&-2&-2        &{}\leq 0& &          4&2&{}   &-1&-2&-1&-2&-3        &{}\leq 1\\  
 2&1&-1&-1&-1&{}    &-2&-2        &{}\leq 0&   &          4&3&-1&-1&-2&-1&-2&-3        &{}\leq 1\\    
 3&2&-1&-1&-1&-1&-2&-3        &{}\leq 0&       &        4&3&-1&{}   &-1&-2&-3&-3        &{}\leq 1\\  
 3&2&-1&{}   &-1&-2&-2&-3        &{}\leq 0&    &         4&3&{}   &-1&-1&-2&-3&-3        &{}\leq 1\\  
 3&2&{}   &-1&-1&-2&-2&-3        &{}\leq 0&    &         4&2&-1&-1&-2&{}    &-3&-3        &{}\leq 1\\  
 3&2&-1&-1&{}    &-2&-3&-3        &{}\leq 0&   &          4&4&-1&-1&-1&-3&-3&-3        &{}\leq 1\\    
 3&3&-1&-1&-1&-2&-3&-3        &{}\leq 0&       &        4&3&-1&-1&-1&-1&-2&-4        &{}\leq 1\\    
 6&4&-1&-1&-2&-3&-4&-6        &{}\leq 0&       &        4&3&-1&{}   &-1&-2&-2&-4        &{}\leq 1\\  
 4&3&-1&-1&-1&-2&-3&-4        &{}\leq 0&       &        4&3&{}   &-1&-1&-2&-2&-4        &{}\leq 1\\  
 5&3&-1&-1&-2&-2&-3&-5        &{}\leq 0&       &        7&5&-1&-1&-2&-3&-4&-6        &{}\leq 1\\    
  &1&-1&{}   &{}    &{}    &{}    &-1        &{}\leq 1&&       5&4&-1&-1&-1&-3&-3&-4        &{}\leq 1\\    
  &1&{}   &-1&{}    &{}    &{}    &-1        &{}\leq 1&&       5&5&-1&-1&-2&-3&-4&-4        &{}\leq 1\\    
  &1&{}   &{}   &-1&{}    &{}    &-1        &{}\leq 1& &      5&4&-1&-1&-1&-2&-3&-5        &{}\leq 1\\    
  &1&{}   &{}   &{}    &-1&-1&{}            &{}\leq 1& &      6&4&-1&-1&-2&-2&-3&-5        &{}\leq 1\\    
  &1&{}   &{}   &{}    &-1&{}    &-1        &{}\leq 1& &      6&4&-1&{}   &-2&-3&-4&-5        &{}\leq 1\\  
  &1&{}   &{}   &{}    &{}    &-1&-1        &{}\leq 1& &      6&4&{}   &-1&-2&-3&-4&-5        &{}\leq 1\\  
 2&1&{}   &{}   &{}    &-1&-1&-1        &{}\leq 1&     &    3&2&-1&{}   &-1&{}    &-1&-2        &{}\leq 2\\
  &1&-1&{}   &-1&{}    &-1&{}            &{}\leq 1&    &     3&2&{}   &-1&-1&{}    &-1&-2        &{}\leq 2\\
  &1&{}   &-1&-1&{}    &-1&{}            &{}\leq 1&    &     3&2&-1&-1&{}    &{}    &-1&-3        &{}\leq 2\\
 2&1&-1&{}   &-1&{}    &-1&-1        &{}\leq 1&        &   3&3&-1&-1&-1&{}    &-1&-3        &{}\leq 2\\  
 2&1&{}   &-1&-1&{}    &-1&-1        &{}\leq 1&        &   4&3&-1&-1&-1&-1&-1&-3        &{}\leq 2\\    
  &1&-1&-1&-1&-1&{}    &{}            &{}\leq 1&       &    5&3&-1&{}   &-2&-1&-2&-4        &{}\leq 2\\  
 2&1&-1&-1&-1&-1&{}    &-1        &{}\leq 1&           &  5&3&{}   &-1&-2&-1&-2&-4        &{}\leq 2\\  
 2&2&-1&-1&-1&-1&-1&-1        &{}\leq 1&               &5&4&-1&-1&-2&-1&-2&-4        &{}\leq 2\\    
 2&1&{}   &{}   &-1&{}    &-1&-2        &{}\leq 1&     &    5&4&-1&{}   &-1&-2&-3&-4        &{}\leq 2\\  
 2&2&-1&{}   &{}    &-1&-1&-2        &{}\leq 1&        &   5&4&{}   &-1&-1&-2&-3&-4        &{}\leq 2\\  
 2&2&{}   &-1&{}    &-1&-1&-2        &{}\leq 1&        &   6&5&-1&-1&-1&-3&-3&-5        &{}\leq 2\\
 2&2&{}   &{}   &-1&-1&-1&-2        &{}\leq 1\\
 \bottomrule
\end{array}$
\end{center}
\end{table}

\begin{table}
\caption{Full description of the subproblems of Example~\ref{ex:branching}}\label{table:ex-branches}
\begin{center}
\small
$\begin{array}{*6{r}@{\;}lc@{\qquad}*6{r}@{\;}l}
\toprule
\multicolumn{7}{c}{\text{Branch $x_7 = 0, x_8 = 0$}} 
&& \multicolumn{7}{c}{\text{Branch $x_7 = 1, x_8 = 0$}}\\
\cmidrule{1-7}\cmidrule{9-15}
c_1 & c_2 & c_3 & c_4 & c_5 & c_6 & \multicolumn{1}{c}{\gamma}
&& c_1 & c_2 & c_3 & c_4 & c_5 & c_6 & \multicolumn{1}{c}{\gamma}\\
\cmidrule{1-7}\cmidrule{9-15}
 1&    &   &   &-1&            &\leq 0       &&  1&    &-1&   &-1&-1        &\leq 0 \\
 1&    &   &   &   &-1        &\leq 0       &&  1&    &   &-1&-1&-1        &\leq 0  \\
 1&    &-1&-1&   &            &\leq 0      &&  1&    &   &   &   &           &\leq 1\\ 
 2&1&   &   &-1&-1        &\leq 0     &&  1&1&-1&   &   &           &\leq 1         \\
 1&1&-1&   &   &-1        &\leq 0     &&  1&1&   &-1&   &           &\leq 1         \\
 1&1&   &-1&   &-1        &\leq 0     &&  1&1&   &   &-1&           &\leq 1         \\
 2&1&-1&-1&-1&            &\leq 0    &&  1&1&   &   &   &-1        &\leq 1          \\
 3&2&-1&-1&   &-2        &\leq 0    &&  3&2&-1&   &-1&-1        &\leq 2             \\
 4&3&-1&-1&-1&-2        &\leq 0   &&  3&2&   &-1&-1&-1        &\leq 2               \\
\midrule
\multicolumn{7}{c}{\text{Branch $ x_7= 0,  x_8= 1$}} 
&& \multicolumn{7}{c}{\text{Branch $ x_7= 1,  x_8= 1$}}\\
\cmidrule{1-7}\cmidrule{9-15}
c_1 & c_2 & c_3 & c_4 & c_5 & c_6 & \multicolumn{1}{c}{\gamma}
&& c_1 & c_2 & c_3 & c_4 & c_5 & c_6 & \multicolumn{1}{c}{\gamma}\\
\cmidrule{1-7}\cmidrule{9-15}
 1&   &   &   &   &           &\leq 1           &&  1&1&-1&-1&-1&-1 &\leq 1 \\
 1&1&   &   &   &-1        &\leq 1         &&\\
 1&1&-1&   &-1&           &\leq 1        &&\\
 1&1&   &-1&-1&           &\leq 1        &&\\
\bottomrule
\end{array}$
\end{center}
\end{table}

This example illustrates why branching is such an important tool in solving
mixed integer programs. The question emerges how to obtain branching decisions
such that the polyhedral description for each of the branches becomes as easy
as possible.  Thus, when we compare branching decisions in our experiments, we 
shall use the following definition.
\begin{definition}
  The \emph{complete description size} of an $n$-way branching decision
  is defined as the sum of the numbers of facets in the complete descriptions
  of the $n$ subproblems.
\end{definition}
Clearly this definition should only be used for comparing branching decisions
with an equal number of subproblems.  For our experiments, we used {\small
  PORTA} \cite{porta}, version 1.3, to enumerate the feasible solutions and to compute the facet
description of their convex hull.  As the computation times for problems of
higher dimension would be prohibitive, we had to restrict ourselves to
experiments with very small integer programs.  Specifically, we generated
dense 0/1 problems with twelve binary variables and two rows.  The four test
instances are shown in Table~\ref{tab:tworowexamples}. 
\begin{table}[tbp]
  \caption{Randomly generated problem instances.  Instances 1~and~2 have
    been generated randomly by drawing the coefficients independently and
    uniformly from the set $\{-20,\dots,+20\}$.  The right-hand side is
    always~0.
    Instance~3 has been modified manually, so that the first three variables have
    identical coefficients.  Finally, instance~4 is a variation of instance~3
    where the coefficients of the first three variables are very close to each
    other.}
  \label{tab:tworowexamples}
  \small
  \begin{center}
    $\begin{array}{*{13}r}
      \toprule
      \multicolumn{12}{c}{\text{Matrix data}}\\
      \cmidrule{1-12}
      A_{1} & A_{2} & A_{3} & A_{4} & A_{5} & A_{6} & A_{7} & A_{8} & A_{9} &
      A_{10} & A_{11} & A_{12} & \ve b \\
      \midrule
      \multicolumn{13}{c}{\text{Instance 1}}\\
      \midrule
      11&-7&9&10&-2&7&14&-15&4&-5&-2&-19&\leq0\\
      6&18&-4&-9&17&-11&5&-12&5&3&-18&7&\leq0\\
      \midrule
      \multicolumn{13}{c}{\text{Instance 2}}\\
      \midrule
      3&-7&0&8&12&-1&7&-14&13&20&-18&2&\leq0\\
      9&11&-13&19&8&-15&-5&3&7&18&-6&-10&\leq0\\
      \midrule
      \multicolumn{13}{c}{\text{Instance 3}}\\
      \midrule
      7&7&7&15&-21&-15&-23&-12&12&-6&11&10&\leq0\\
      10&10&10&-21&4&-3&4&13&-1&-14&2&-6&\leq0\\
      \midrule
      \multicolumn{13}{c}{\text{Instance 4}}\\
      \midrule
      7&6&7&15&-21&-15&-23&-12&12&-6&11&10&\leq0\\
      10&10&9&-21&4&-3&4&13&-1&-14&2&-6&\leq0\\
      \bottomrule
    \end{array}$
  \end{center}
\end{table}

\paragraph{On the limitations of current LP-based branching schemes.}

A single-variable branching scheme, which is used in today's branch-and-cut
systems, is usually driven by information obtained from the current LP
relaxation (``most infeasible variable selection''), by lookahead-based
techniques (``strong branching''), and history-based prediction
(``pseudo-cost branching'').  There is a large class of problems
that are extremely difficult to solve for current branch-and-cut systems
because none of the above criteria provides a meaningful basis for a
branching decision.  An extreme example for this are the market split
instances by  
Cornu\'ejols and Dawande \cite{corn-dawande:98}:  Here the LP relaxations of
all subproblems have the value~$0$, until most of the variables have already been
fixed.  However, it was shown that branch-and-bound is indeed the right tool
for solving the market split instances:  While LP-based single-variable
branching fails, it is very successful to branch on certain
general disjunctions that can be derived from the problem
structure via lattice basis reduction~\cite{aardal-hurkens-lenstra-2000}.
Though this technique has proved very successful for solving market split
problems~\cite{ABHLS2000} and also for the so-called
banker's problem~\cite{louveaux-wolsey-2002:basis-reduction}, it has not
become a general tool for branch-and-cut systems.  

We also refer to the recent
work~\cite{karamanov-cornuejols-2005:branching-general} where a branching
method along general disjunctions is proposed.  Here the quality of a
disjunction (branching direction) is measured by the depth of the intersection
cut corresponding to the disjunction; among the best disjunctions, strong
branching is used to select one.  The computational results for many benchmark
problems from MIPLIB are very promising.  However, for a few instances the
proposed branching scheme fails to close any gap.  This includes the market split
instances \texttt{markshare1} and \texttt{markshare2}. 

\paragraph{A branching scheme based on value disjunctions.}
We propose a new branching scheme based on value disjunctions, which we hope
is general enough to be useful as a branching scheme for general integer
programs.  It is purely based on the analysis of the structure of the integer
program, and is designed to complement the above mentioned LP-based prediction
methods.  

The basic idea of the new branching scheme is to partition the set~$N$ of
problem variables into blocks $N_i$ and to move over to the extended formulation given
by the value disjunction.  In addition to the original variables, we can then
branch on the newly introduced binary variables.  In fact, because exactly one
binary variable of each block can be set to~1, we can perform {\small SOS}
branching on these variables.  The question, of course, is how to construct a
suitable partition of~$N$.  
\begin{claim}
  One should choose a set of variables
  whose columns are structurally similar and perform a value disjunction
  according to a relaxation where we replace the original coefficients by simpler
  ones.
\end{claim}

For our experiment, we decided to pick three of the twelve binary variables,
$x_i$, $x_j$, $x_k$, say.  We then add the (redundant) constraint $x_i + x_j +
x_k\leq 3$.  When we construct a value disjunction with respect to this
constraint, we need to introduce four variables $y_0$, $y_1$, $y_2$, $y_3$,
corresponding to the possible values of the form $x_i + x_j + x_k$.
Performing {\small SOS} branching on $y_0+y_1+y_2+y_3=1$ yields four
subproblems.  To compute the complete description size of the value
disjunction branching on $x_i$, $x_j$, $x_k$, we sum up the numbers of facets in
each of these four subproblems.  To make a comparison with traditional
single-variable branching, we need to consider a branching strategy that
yields the same number of subproblems.  To this end, we pick two
original variables, $x_i$, $x_j$ say, and consider the subproblems where we fix
these variables to the possible values. 

We next defined a ``ranking formula'' for the selection of the three variables
$x_i$, $x_j$, $x_k$ that give rise to the value disjunction.  Let $A_i$,
$A_j$, $A_k$ denote the columns of these 
variables.  Then let
\begin{equation*}
  R(\{i,j,k\}) = \min\nolimits_{r=1}^2 \frac{{\left(\max\{A_{r,i}, A_{r,j},
        A_{r,k}\} - \min\{A_{r,i}, A_{r,j}, A_{r,k}\}\right)}^2 } 
  {2 + \left|\mathop{\mathrm{med}}\{A_{r,i}, A_{r,j}, A_{r,k}\}\right|}
\end{equation*}
where $\mathrm{med}\{A_{r,i}, A_{r,j}, A_{r,k}\}$ denotes the median of the
three values.  The formula was designed so that (i)~columns that have
``similar'' coefficients in at least one of the rows yield a low (good)
result; (ii)~columns with large coefficients yield a low result.
The rationale of this ranking is that, intuitively, the value disjunction for
a selection of similar columns should lead to simpler subproblems; also
columns with large coefficients should have a larger impact on the rest of the
problem than columns with small coefficients.
\begin{example}
  For test instance~4, selecting the variables $x_1$, $x_2$, $x_3$ has the rank
  $R(\{1,2,3\})= 0.083$; selecting the variables $x_7$, $x_9$, $x_{10}$ has
  the rank $R(\{7,9,10\})=108$. 
\end{example}

\begin{figure}
  \pgfdeclareimage[width=.84\textwidth]{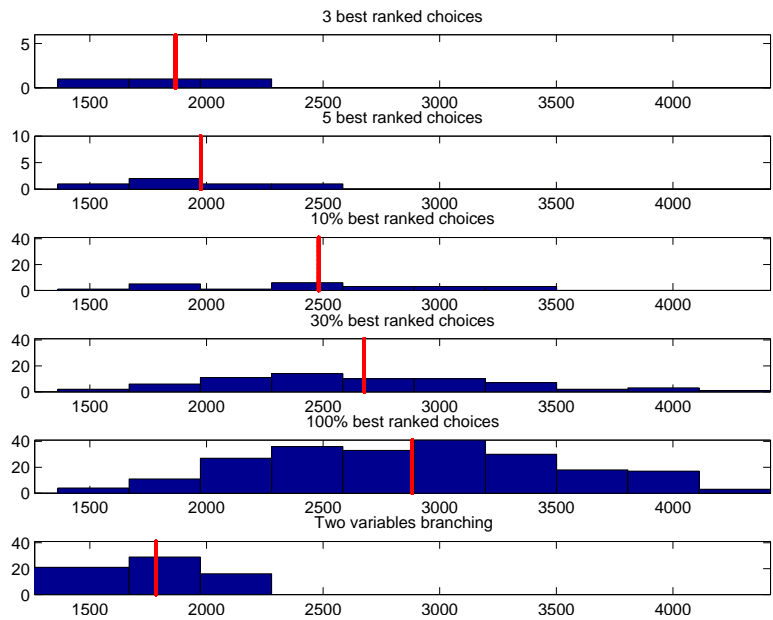}{knapsack1-histograms}
  \begin{center} 
    \pgfuseimage{knapsack1-histograms} 
  \end{center} 
  \caption{Branching on value disjunctions vs.\ 2-variable branching
    (instance 1).  The figure shows histograms of the total number of facets
    in the subproblems; the vertical line is the average.}\label{knapsack1-histograms}
\end{figure}

\begin{figure}
  \pgfdeclareimage[width=.84\textwidth]{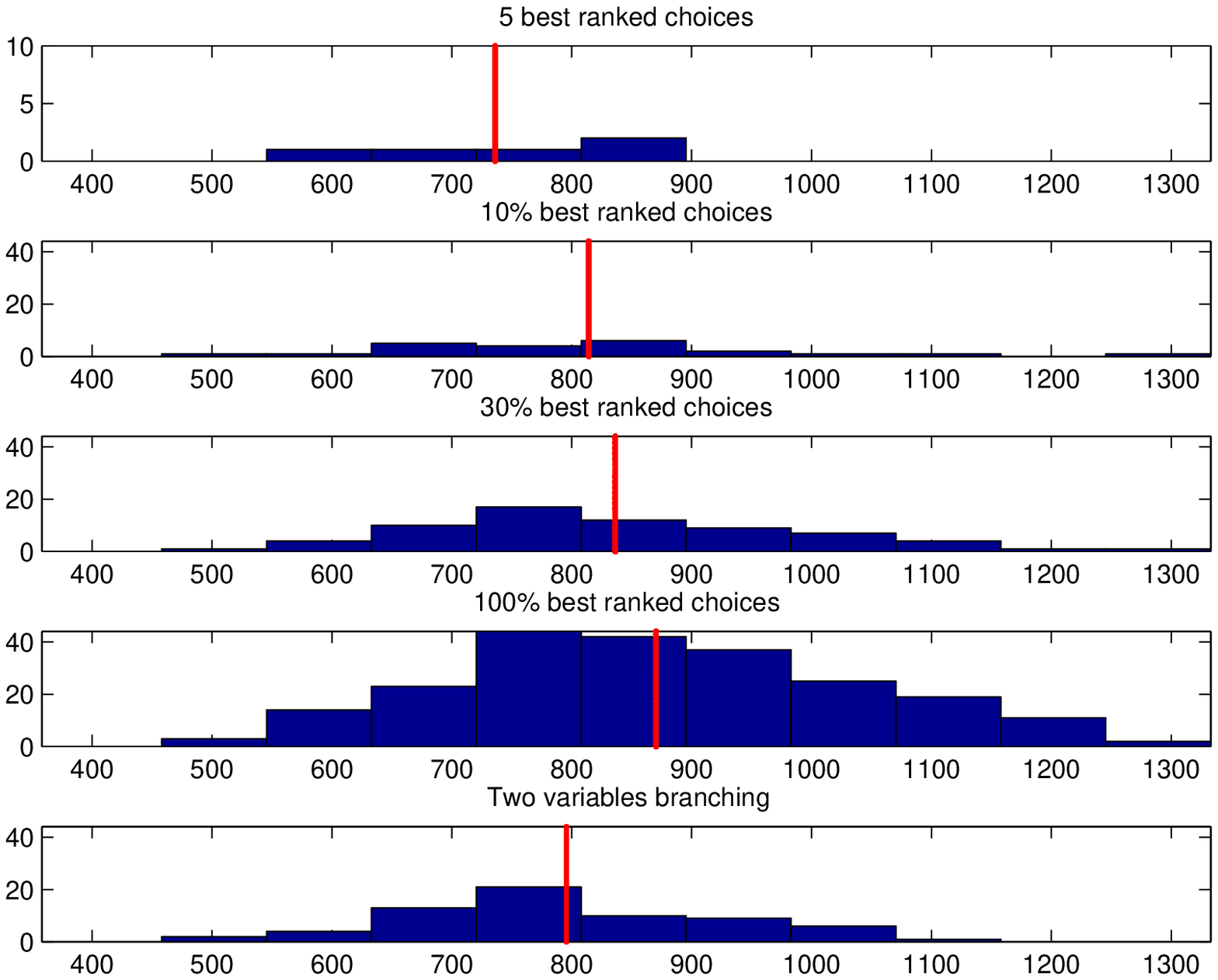}{instance2-histograms}
  \begin{center} 
    \pgfuseimage{instance2-histograms}
  \end{center}
  \caption{Branching on value disjunctions vs.\ 2-variable branching (instance
    2)} 
  \label{instance2-histograms}
\end{figure}

\begin{figure}
  \pgfdeclareimage[width=.84\textwidth]{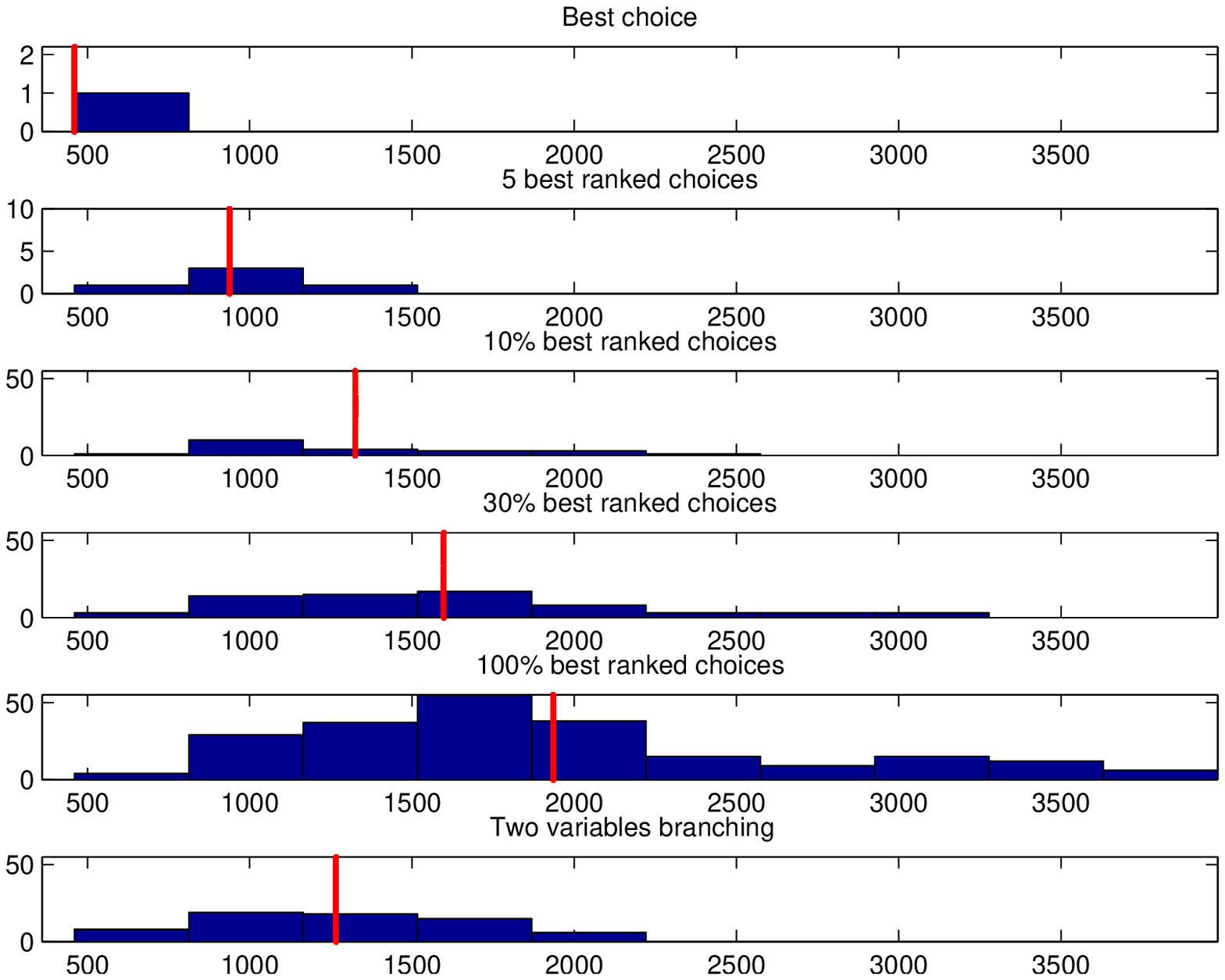}{structured-instance2-histograms}
  \begin{center} 
  \pgfuseimage{structured-instance2-histograms}
\end{center}
  \caption{Branching on value disjunctions vs.\ 2-variable branching (instance 3)} 
  \label{structured-instance2-histograms}
\end{figure}

\begin{figure}
  \pgfdeclareimage[width=.84\textwidth]{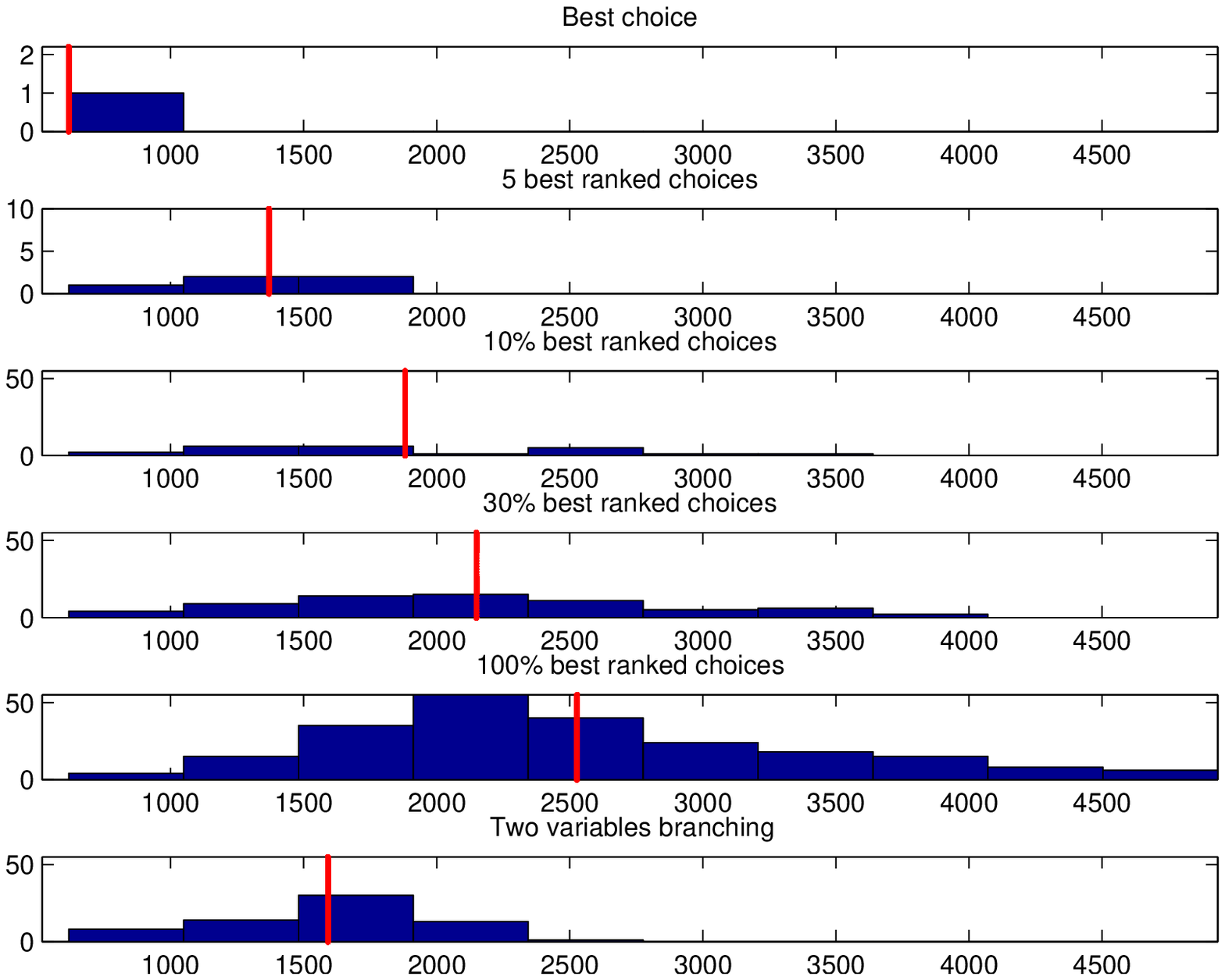}{structured-instance3-histograms}
  \begin{center} 
    \pgfuseimage{structured-instance3-histograms}
  \end{center}
  \caption{Branching on value disjunctions vs.\ 2-variable branching (instance 4)} 
  \label{structured-instance3-histograms}
\end{figure}

For all possible branching decisions (i.e., the $\binom{12}{3}$ choices of three variables), we
now computed the rank and the complete description size.  We grouped the
branching decisions according to their rank into sets of the 5 best ranked, 10~\% best
ranked, 30~\% best ranked, etc.\ choices.  For each of the test instances,
we show histograms of the complete description sizes corresponding to branching
decisions within these rankings in
Figures~\ref{knapsack1-histograms}--\ref{structured-instance3-histograms}.   
As a comparison, the bottom part in each figure shows a histogram of the
complete description sizes obtained by the $\binom{12}{2}$ possible choices
for two-variable branching.
In each histogram the vertical line shows the average (arithmetic mean) of the complete
description sizes.\smallbreak

\relax From the computational results, we can draw the following conclusions:
\begin{enumerate}
\item It is possible to use the rank formula to predict which branching
  decisions will lead to low complete description sizes.
\item For instances 1~and~2 that do not contain selections of very low rank,
  two-variable branching performs better than branching of value disjunctions.
  However,  instances 3~and~4 that contain selections of very low rank, it is possible
  to take branching decisions that are better than two-variable branching
  decisions by making use of the rank formula.
\end{enumerate} 

We have to remark that there is room for improvement of the proposed
ranking formula.  Clearly it needs to be generalized for blocks of different
cardinalities.  It would also need adjustment for unequally scaled
rows.  

\paragraph{Value disjunction branching on larger problems.}

Based on the evidence obtained with the above experiments, we tried to use the
new branching scheme to solve larger test problems.  Our set of test instances
consists of instances with several dense rows (multi-knapsack problems).  We
focused on problems where the solutions to LP relaxations of subproblem only
give little information for taking branching decisions.  The test instances
are: 
\begin{itemize}
\item Six randomly generated market split instances with 35 and 40 variables.
\item The models \texttt{mas74} and \texttt{mas76} from
  the MIPLIB.
\end{itemize}

It seems difficult to apply Theorem  
\ref{structure theorem for value disjunction} directly to these problems. The reason is that typically many constraints in a model are present. In this case the probability that we can come up with a block decomposition such that some values repeat, is quite low. Hence, one may expect that in such cases the value-reformulation requires to introduce as many variables as we have subsets in each of the elements of the partition $N_1,\ldots, N_K$.
Therefore, we decided to perform the following steps: 
\begin{enumerate} 
\item We consider one of the dense rows at a time. We add a relaxation of this
  row that we obtain by replacing the coefficients by simpler ones. From the row 
  \begin{equation*} 
    \sum_{i=1}^n a_i x_i + \sum_{j=1}^d g_j w_j \leq  b, 
  \end{equation*} 
  we generate the relaxation 
  \begin{equation*} 
    \sum_{i=1}^n f(a_i) x_i \leq M, 
  \end{equation*} 
  where $f(x)$ is a non-linear function of the type 
  \begin{equation*} 
    f(x) =  
    \left\{ 
      \begin{array}{ll} 
        1 & \text{ if } x \geq U\\ 
        0 & \text{ if } L < x < U\\ 
        -1 & \text{ if } x \leq L. 
      \end{array} 
    \right. 
  \end{equation*}   
\item We reformulate the problem using a value disjunction for each of the new
  rows separately.
\item Finally, we manually perform {\small SOS} branching on the new
  variables.  Then we solve each of the subproblems with the standard branch-and-cut
  system {\small CPLEX}~9.1~\cite{ilog-cplex-uuh} using the default settings of the Callable Library.  We use
  the optimal solution value from a subproblem as a 
  primal bound for the remaining subproblems.
\end{enumerate} 
The results of this approach on the set of test instances are shown in
Table~\ref{cornies}.  It can be seen that the approach provides a clear gain
on all these instances.  Both the number of nodes and the computation times
are reduced in comparison to the performance of {\small CPLEX}~9.1 (with the default
settings of the Callable Library) on the original problem.
 
\begin{table}[tbp] 
  \caption{Branching on value disjunctions for the market split and
    \texttt{mas} instances. 
    Computation times are given in CPU seconds on a Sun Fire V890 with 1200\,MHz UltraSPARC-IV processors} 
  \label{cornies} 
  \small
  \begin{center} 
\def~{\hphantom0} 
\begin{tabular}{lcccccc}
\toprule 
 &  &  & \multicolumn{2}{c}{CPLEX 9.1} & 
\multicolumn{2}{c}{Value Disjunctions}\\ 
\cmidrule(lr){4-5}\cmidrule(lr){6-7}
\multicolumn{1}{c}{Name} & Rows & Cols & Nodes ($10^6$)& Time (s) & Nodes ($10^6$)& Time (s)\\ 
\midrule 
\texttt{corn535-1}&~5 & ~40 & ~13.8\hphantom{00} & ~2\,431 & ~~3.8\hphantom{00} & ~~809 \\ 
\texttt{corn535-2}&~5 & ~40 & ~11.9\hphantom{00} & ~2\,084 & ~~4.2\hphantom{00} & ~~865 \\ 
\texttt{corn535-3}&~5 & ~40 & ~17\hphantom{.000} & ~2\,946 & ~~9.8\hphantom{00} & ~1\,970 \\ 
\texttt{corn540-4}&~5 & ~45 & 321\hphantom{.000} & 55\,918 & 105\hphantom{.000} & 20\,873 \\ 
\texttt{corn540-5}&~5 & ~45 & 231\hphantom{.000} & 39\,787 & ~87\hphantom{.000} & 17\,267 \\ 
\texttt{corn540-6}&~5 & ~45 & 188\hphantom{.000} & 30\,532 & ~97\hphantom{.000} & 19\,162 \\ 
\texttt{mas74}   & 13 & 151 & ~~4.4\hphantom{00} & ~2\,463 & ~~1.2\hphantom{00} & ~1\,194 \\ 
\texttt{mas76}   & 12 & 151 & ~~0.667            & ~~\,289 & ~~0.063     & ~~\,~35 \\ 
\bottomrule
\end{tabular} 
\end{center} 
\end{table} 
 
\clearpage
\bibliographystyle{amsplain} 
\bibliography{weismantel,iba-bib} 
 
\end{document}